\documentclass[11pt,twoside]{article}
\usepackage{latexsym,amssymb,amsmath}
\usepackage{amsfonts}
\usepackage{longtable}
\RequirePackage{graphics,epsfig,psfrag}
\def\abs#1{\left \vert #1 \right \vert}

\def\RR{{\bf R}} 
\def\ZZ{{\bf Z}} 
\def\NN{{\bf N}} 
\def\Mod#1{\,(\hbox{\rm mod}\,#1)}
\def\cn{\hbox{\rm cn}\,}
\def\H{{\rm H}}

\def\phi{\varphi}
\def\eps{\varepsilon}

\def\cC{{\cal C}}

\def\cT{{\cal T}}
\def\pn{\medskip\par\noindent}
\def\bi{\vspace{-2pt}\begin{itemize}\itemsep -2pt plus 1pt minus 1pt}
\def\ei{\end{itemize}\vspace{-4pt}}
\def\bn{\vspace{-2pt}\begin{enumerate}\itemsep -2pt plus 1pt minus 1pt}
\def\en{\end{enumerate}\vspace{-4pt}}
\newcommand{\Pf}{{\em Proof}. }

\newcommand{\EPf}{\hbox{}\hfill$\Box$\vspace{.5cm}}
\def\[#1\]{\begin{eqnarray}#1\end{eqnarray}}
\def\$#1\${\begin{eqnarray*}#1\end{eqnarray*}}

\def\pent#1#2{\lfloor\frac{#1}{#2}\rfloor}

\def\sign#1{{\rm sign}\,\bigl( #1 \bigr)}

\def\abs#1{\left \vert #1 \right \vert}

\def\frac#1#2{{\textstyle{{#1} \overwithdelims.. {#2}}}}
\def\Frac#1#2{{\displaystyle{{#1} \overwithdelims.. {#2}}}}
\def\p{{\small $+$}}
\def\m{{\small $-$}}
\catcode`\@=11
\def\system#1{\left\{\null\,\vcenter{\openup\jot\m@th
\ialign{
\strut\hfil$\displaystyle{##}$&
$\,\displaystyle{{}##}\,$\hfil&&
\strut\hfil$\,\displaystyle{##}$&
$\,\displaystyle{{}##}\,$\hfill
\hfil\crcr#1\crcr}}\right.}

\def\cmatrix#1{\left [
\null\,\vcenter{
\ialign{
\hfil${##}\ $\hfil &
\hfil$\ {##}\ $\hfil&&
\hfil$\ {##}\ $\hfil&
\hfil$\ {##}$\hfil
\crcr#1\crcr}}\right ]}

\def\@opargbegintheorem#1#2#3{\par\addvspace{6pt plus3pt minus2pt}%
    \def\@tempa{#3}%
    \noindent{\bf #1 #2 \ifx\@tempa\empty\unskip\else\unskip\ (#3).\fi\hskip.5em}\csname#1font\endcsname\ignorespaces
\ignorespaces}

\def\@endtheorem{\par\addvspace{6pt plus3pt minus2pt}}
\def\@begintheorem#1#2#3{\par\addvspace{8pt plus3pt minus2pt}%
              \noindent{\csname#1headfont\endcsname#1\ \ignorespaces#3 #2.}%
              \csname#1font\endcsname\hskip6pt\ignorespaces}
\def\@endtheorem{\par\addvspace{8pt plus3pt minus2pt}\@endparenv}

\catcode`\@=12

\newtheorem{theorem}{Theorem}[section]
\newtheorem{thm*}{Theorem}
\newtheorem{corollary}[theorem]{Corollary}
\newtheorem{lemma}[theorem]{Lemma}
\newtheorem{proposition}[theorem]{Proposition}
\newtheorem{definition}[theorem]{Definition}
\newtheorem{remark}[theorem]{Remark}
\newtheorem{example}[theorem]{Example}

\date{\today}
\begin{document}
\pagestyle{myheadings} \markboth{P. -V. Koseleff, D. Pecker}{{\em
Harmonic knots}}
\title{ Harmonic knots}
\author{P. -V. Koseleff, D. Pecker\medskip}
\maketitle
\begin{abstract}
The harmonic knot $\H(a,b,c)$
is parametrized as
 $K(t)= (T_a(t) ,T_b (t), T_c (t))$ where
$a$, $b$ and $c$ are  pairwise coprime integers and $T_n$ is the degree $n$
Chebyshev polynomial of the first kind.
We classify the harmonic knots $\H(a,b,c)$ for $ a \le 4. $
We study the knots $\H (2n-1, 2n, 2n+1),$
the knots $\H(5,n,n+1),$ and give a table
of the simplest harmonic knots.
\pn {\bf keywords:}{
Long knots, polynomial curves, Chebyshev curves, rational knots, continued fractions}
\end{abstract}
\begin{center}
\parbox{12cm}{\small
\tableofcontents
}
\end{center}
\vspace{1cm}
\section{Introduction}
A harmonic curve (or Chebyshev curve) is defined to be
a curve which admits a parametrization
$x=T_a(t), \, y= T_b(t), \, z=T_c(t)$ where $t \in \RR  ,$  $a$, $b$ and $c$ are
 integers,   and $T_n(t)$ are  the  Chebyshev polynomials
defined by $T_n(\cos t) = \cos nt.$
A harmonic knot is a nonsingular harmonic curve, it is a long knot.
In 1897  Comstock proved that a harmonic curve is a knot if and only
if $a,b,c$ are pairwise  coprime integers (\cite{Com,KP3,Fr}).
\pn
We observed in \cite{KP1} that  the trefoil can be parametrized by
Chebyshev polynomials: $x=T_3(t),\, y=T_4(t), \,z= T_5(t)$.
This led us to study harmonic knots in \cite{KP3}.
\pn
Harmonic knots are polynomial analogues of the famous Lissajous knots
(\cite{BDHZ,BHJS,Cr,HZ,JP,La1,La2}).
However, they are very different:
there are  only two known examples of knots which are both Lissajous and harmonic, the knots $ 5_2$ and $7_5.$
\pn
We proved in \cite{KP3} that the harmonic knot $\H(a,b, ab-a-b)$ is alternating, and deduced
that there  are infinitely many amphicheiral harmonic knots and
infinitely many strongly invertible harmonic knots.
We also proved  that the torus knot $T(2, 2n+1)$ is the
harmonic knot $ \H(3,3n+2,3n+1)$.
\pn
The harmonic knots $\H(3,b,c)$ are classified
in \cite{KP4};
they  are two-bridge knots and their Schubert fractions
$\Frac\alpha\beta$ satisfy $\beta^2 \equiv \pm 1 \Mod\alpha$.
\pn
In this article, we give the classification of the harmonic knots
$\H(4,b,c)$ for $b$ and $c$ coprime odd integers. We also study some
infinite families of harmonic knots for $a\ge 5.$
\pn
In section {\bf \ref{cf}} we recall the Conway notation for  two-bridge knots,
and the computation of their Schubert fractions.
The knots $ \H (4,b,c)$ are two-bridge knots, and
their Schubert fractions are given by  continued fractions
 of the form  $[\pm 1,\pm 2,  \ldots, \pm 1, \pm 2].$
In section  {\bf \ref{harmonic}} we compute the Schubert fractions of
the  knots $\H(4,b,c),$ and we deduce their classification.
\pn
{\bf Theorem \ref{h4bc}.}
{\em
Let $b$ and $c$ be relatively prime odd integers, and let $K= \H(4,b,c).$
There is a unique pair
$(b',c')$
such that (up to mirroring)
$$K= \H(4,b',c'), \    b'<c'<3b', \  b' \not\equiv c' \Mod 4.$$
$K$ has a   Schubert fraction
$\Frac \alpha \beta $  such that
$ \beta^2 \equiv \pm 2 \Mod \alpha.$
Furthermore, there is an algorithm to find $(b', c'),$
and the crossing number of $K$ is $ N= (3b'+c'-2)/4. $
}
\pn
We notice that the trefoil is the only knot which is both of form
$\H(3,b,c)$ and $\H(4,b,c).$
In section {\bf \ref{special}}  we study  some families of harmonic
knots $\H (a,b,c)$ with $a \ge 5.$  In general their bridge number
 is greater than two, this is why the following result is
surprising.
\pn
{\bf Theorem \ref{surprise}.}
{\em
The harmonic knot $ \H(2n-1,2n,2n+1)$ is isotopic to the two-bridge harmonic knot
$\H (4, 2n-1,2n+1),$ up to mirror symmetry.
}
\pn
We also find an infinite family of two-bridge harmonic knots which are
not of the form $ \H (a,b,c)$ for $a \le 4$:
\pn
{\bf Theorem \ref{h5}.}\\
{\em
The knot $ \H ( 5, 5n+1, 5n+2)$ is the two-bridge knot of Conway form
$C ( 2n+1, 2n).$\\
The knot $ \H ( 5, 5n+3, 5n+4)$ is the two-bridge knot of Conway form
$ C (2n+1, 2n+2).$\\
Except for $\H(5,6,7)= \H (4,5, 7)$ and $ \H(5,3,4),$
these knots are not of the form $ \H (a,b,c)$ with $a \le 4.$}
\pn
Then, we identify  the  knots  $ \H( a,b,c) $ for $  (a-1)(b-1) \le 30.$ Our examples
show that harmonic knots are not necessarily prime, nor reversible.
\section{Continued fractions and two-bridge knots}\label{cf}
\begin{figure}[th]
\begin{center}
{\scalebox{.7}{\includegraphics{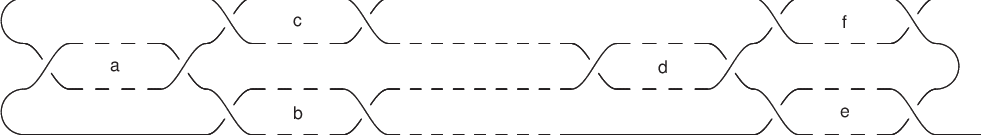}}}
\end{center}
\vspace{-12pt}
\caption{A knot isotopic to $C(b_1,a_1+c_1,b_2,a_2+c_2,\ldots, b_n,a_n+c_n)$}
\label{conway4}
\end{figure}

A two-bridge knot (or link) admits a diagram in Conway form.
This form, denoted by
$C(a_1, a_2, \ldots, a_n)$  where $a_i \in \ZZ$, is explained by
the following picture (see \cite{Con}, \cite[p. 187]{Mu}).
\psfrag{a}{\small $a_1$}\psfrag{b}{\small $a_2$}%
\psfrag{c}{\small $a_{n-1}$}\psfrag{d}{\small $a_{n}$}%
\begin{figure}[th]
\begin{center}
{\scalebox{.8}{\includegraphics{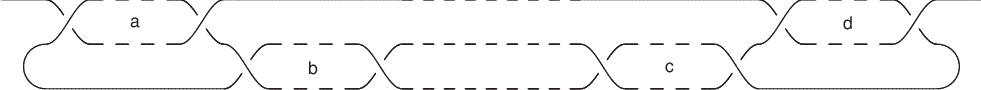}}}\\[30pt]
{\scalebox{.8}{\includegraphics{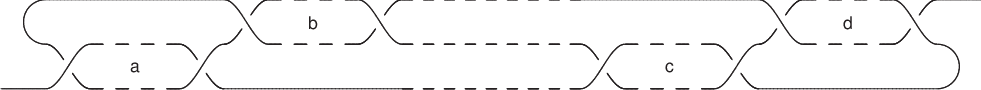}}}
\end{center}
\vspace{-12pt}\caption{Conway  forms  for polynomial knots ($n$
  odd, and $n$ even)}
\label{conways3}
\end{figure}
The number of twists is denoted by the integer
$\abs{a_i}$, and the sign  of $a_i$ , called the Kauffman sign, is defined
as follows: if $i$ is odd, then the right twist is positive,
if $i$ is even, then the right twist is negative.
In Figure \ref{conways3} the $a_i$ are positive (the $a_1$ first twists are right twists).
\pn
The two-bridge knots (or links) are classified by their Schubert fractions
$$
 \Frac {\alpha}{\beta} =
a_1 + \Frac{1} {a_2 +  \Frac{1} {\cdots +\Frac 1{a_n}}}=
[ a_1, \ldots, a_n], \quad \alpha >0, \   \   a_i \in \ZZ \cup \{ \infty  \} .
$$
    A two-bridge knot (or link) with
Schubert fraction $ \Frac {\alpha}{\beta}$ is denoted by $S \bigl( \Frac {\alpha}{\beta} \bigr).$
The two-bridge  knots (or links)
$ S (\Frac {\alpha} {\beta} )$ and $ S( \Frac {\alpha ' }{\beta '} )$ are equivalent
if and only if $ \alpha = \alpha' $ and $ \beta' \equiv \beta ^{\pm 1} ( {\rm mod}  \  \alpha).$
If $K= S (\Frac {\alpha}{\beta} ),$ its mirror image is
$ \overline{K}= S ( \Frac {\alpha}{- \beta} ).$
\pn
We shall   study knots with a diagram illustrated by Figure \ref{conway4}.
\psfrag{a}{\small $b_1$}\psfrag{b}{\small $a_1$}\psfrag{c}{\small $c_1$}%
\psfrag{d}{\small $b_n$}\psfrag{e}{\small $a_n$}\psfrag{f}{\small $c_n$}%
In this case, the $a_i$ and the $c_i$ are positive if they are left twists,
the $b_i$ are positive if they are right twists (in our figure  $a_i, b_i, c_i $ are positive).
Such a knot is equivalent to a knot of Conway form
$C(b_1, a_1+c_1, b_2, a_2 +c_2, \ldots, b_n, a_n+c_n )$
(see \cite[ p.~183-184]{Mu}).
\medskip
Our knots have a Chebyshev diagram, that is a (singular) plane Chebyshev curve
$\cC(4,k): x= T_4(t), \, y= T_ k(t)$, and the over/under information at each crossing.
In this case we obtain diagrams of the form illustrated by Figure \ref{conway4}.
Then, by symmetry such a knot  has a Schubert fraction of the form
$[ b_1, d_1, b_2, d_2,  \ldots, b_n, d_n ]$
with $b_i=\pm 1$,  and $d_i = \pm 2$.

\begin{figure}[th]
\begin{center}
\begin{tabular}{c}
\psfrag{n1}{\m}\psfrag{n2}{\m}\psfrag{n3}{\m}%
\psfrag{n4}{\p}\psfrag{n5}{\p}\psfrag{n6}{\p}%
{\scalebox{.9}{\includegraphics{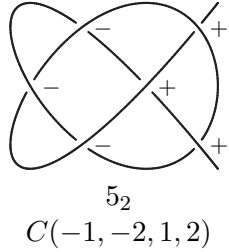}}}\\
$5_2$\\
$C(-1,-2,1,2)$
\end{tabular}
\end{center}
\vspace{-12pt}
\caption{ A Chebyshev diagram  of the harmonic knot $H(4,5,7) = 5_2 $ }\label{diag4}
\end{figure}
\subsection{Continued fractions}
Let $K$ be a two-bridge knot defined by its Conway  form
$C(q_1, q_2, \ldots, q_n ) , $
where   $q_i \in \ZZ$.
It is often possible to obtain directly the crossing number of $K .$
\begin{definition}
Let $r >0$ be a rational number, and $r=[q_1, \ldots, q_n ] $ be a continued fraction
with $q_i \in \NN .$
The {\em crossing number } of $r$ is defined by
$\cn(r)= q_1 + \cdots + q_n. $
\end{definition}
The following result is proved in \cite{KP4}.
\begin{proposition}\label{bireg}
Let $\Frac \alpha\beta = [a_1, \ldots, a_n], \  a_i \in \ZZ$
be a  continued fraction such that
$ a_1a_2>0, \  a_{n-1}a_n >0, $ and without any  two consecutive sign changes in the
sequence $ a_1, a_2, \ldots , a_n.$
Then its crossing number is
\[
\cn(\Frac \alpha\beta) = \sum_{k=1}^n \abs{a_i} - \sharp \{j,
a_ja_{j+1}<0\}.
\label{biregf}
\]
\end{proposition}
\subsection{Continued fractions ${{[\pm 1,\pm 2,\ldots,\pm 1,\pm 2]}}$}
We begin with a useful lemma:
\begin{lemma} \label{unic}
 Let $r= [ 1, 2e_2,e_3, 2e_4, \ldots , e_{2m-1}, 2 e_{2m} ],  \  e_i=\pm 1.$
We suppose that there are no three consecutive sign changes
in the sequence $ e_1, \ldots , e_{2m}.$ Then $ r>0$, and $r>1$ if and only if $e_2=1.$
Here, we use the convention that  $ \infty$ is greater than all rational numbers.
\end{lemma}
\Pf
 By induction on $m.$\\
If $m=1,$ then $r=[1,2]= \Frac 32  \ $ or $ \ r= [1, -2] = \Frac 12 ,$
and the result is true.
\pn
Suppose the result true for $m-1,$ and let us prove it for $m.$
\bn
\item[] First, let us suppose $r= [1,2, 1 \ldots , 2 e_{2m} ]$.
Then $r= [1, 2, y] =  \Frac {3y+1}{2y+1},$ where $y= [1, \pm 2, \ldots ] .$
By induction we have $y>0$ (or $y = \infty$ ), and then $ r > 1.$
\pn
\item[]
Now, let us suppose $ r = [ 1,2,-1,2, \ldots ].$
  If $m=2$, then $r=\infty$ and the result is true.
If $ m \ge 3,$ then $e_5=1$ and
$r= [1,2,-1,2, y]= y+2$ with $y=[1 , \pm 2, \ldots ].$ We have $y>0$ $  ( {\rm or} \  \   y=\infty)$
by induction, and then $r>2>1$ (or $r=\infty$).
\pn
\item[]
Then let us suppose   $r= [1,2, -y]= \Frac {3y-1}{2y-1}= \Frac 32 + \Frac {1}{2(2y-1)}$
with $y= [1,2, \pm1, \ldots].$
Then we have $y>1$ (or $y=\infty$) by induction, and then $ r \ge \frac 32 >1.$
\pn
\item[] Finally, let us suppose $r= [1, -2, \ldots ] .$

\item[]
If $r= [1, -2, -1 , \ldots]$, then $r= [1,-2,-y] = \Frac {y+1}{2y+1} ,$ with
$y= [1, \pm 2 , \ldots ] .$  By induction, we have $y>0$ (or $y= \infty$), and then $ 0<r<1.$
\item[]
If $r= [1, -2, 1 , \ldots]$,
then $r = [ 1 , -2, y]= \Frac{y-1}{2y-1}$ where $ y= [1,2, \pm1, \ldots ].$
By induction we have $y>1$ (or $y=\infty$) and then $ 0<r<1.$
\en
This completes the proof.\EPf

\begin{remark}\label{three}
Because of the identities
$ x= [1,-2,1,-2,x]$ and $ x= [2,-1,2,-1, x],$
we see that the condition on the sign changes is necessary.
\end{remark}
\begin{theorem} \label{th:cf1212}
Let $r=\Frac {\alpha}{\beta} >0 $ be a fraction with $\alpha$ odd and $\beta$ even.
There is a unique
continued fraction expansion
$ r= [1, \pm2, \ldots , \pm1, \pm2 ] $ without three consecutive sign changes.
\end{theorem}
\Pf
The existence of this continued fraction expansion is given in \cite{KPR}.
Its uniqueness is a direct consequence of  Lemma \ref{unic}.
\EPf
\pn
The next result will be useful to describe the continued fractions of
harmonic knots $\H(4,b,c).$
\begin{proposition} \label{palin4}
Let $r= \Frac \alpha \beta$ be a rational number given by a
continued fraction of the form
$ r= [e_1,  2e_2, e_3,  2e_4, \ldots e_{2m-1}, 2 e_{2m} ], \, e_1=1, \,    e_i= \pm 1.$
We  suppose that the sequence of sign changes is palindromic,
that is  $ \  e_k e_{k+1} = e_{2m-k} e_{2m-k+1} \  $ for  $k=1, \ldots, 2m-1.$
Then we have $\beta^2 \equiv \pm 2 \Mod \alpha.$
\end{proposition}
\Pf
We shall use the M\"{o}bius transformations
$$ A(x)=[1,x] =\Frac{x+1}{x+0} , \quad
B(x)=[2,x]=\Frac{2x+1}{x+0}, \quad  S(x)=-x$$
and their matrix notations
$$
A = \cmatrix{1&1\cr1&0}, \
B = \cmatrix{2&1\cr1&0}, \
S = \cmatrix{1&0\cr0&-1},  \
AB= \cmatrix{3&1\cr2&1},  \
ASB= \cmatrix{1&1\cr2&1} .
$$
We shall consider the mapping (analogous to matrix transposition)
$$\tau: \cmatrix {a & b \cr c & d }\mapsto \cmatrix { a & \Frac c2 \cr
  2b & d }.$$
We have $ \tau(XY)= \tau (Y) \tau (X), \ \tau(AB)=AB,   \  \tau(ASB)= ASB$
and $\tau(S) = S $.
\pn
Let $G$ be the M\"{o}bius transformation defined by
$G(z) =  [1,2e_2, e_3,  2e_4, \ldots e_{2m-1}, 2 e_{2m} , z].$
We have $\Frac\alpha\beta = G (\infty).$
Let us write $G=X_1 \cdots X_n$  where $X_i= A, B $ or $S$,
$X_1= A$ and $X_n=B$.
One can suppose that $G$ contains no subsequence of the form $AA, ASA,
BB, SS$ and $BSB$.
Moreover, the palindromic condition means that if $X_i= S,$
then $X_{n+1-i} = S.$
\pn
Let us show that if $P= X_1 \cdots X_n$ is  a product of terms $A, B, S$
having these properties, then $ \tau (P) = P,$
by induction on $s = \sharp\{i, X_i= S\}.$
\pn
If $s=0$ then $P = (AB)^m,$ and $ \tau (P) = \bigl(\tau (AB) \bigr)^m
= (AB)^m= P$.
\pn
Let $k = \min \{i, X_i=S\}.$ Since $X_1 \neq S,$ we have $k \ne 1.$
\bn
\item[]
If $k=2q+1$ then $q \ge 1$ and $ P= (AB)^q S \, P' \, S (AB)^q.$
By induction we have $ \tau ( P')= P', $ and then
$ \tau (P) = \tau \bigl( (AB)^q \bigr) \tau (S) \tau (P') \tau (S)
 \tau \bigl( (AB)^q \bigr) = P .$
\item[]
If $k=2q$ then $P = ( AB)^{q-1} \, (ASB) \, P' \, (ASB) \, (AB)^{q-1}.$
By induction we have $ \tau (P')= P',$ and then $ \tau (P) = P.$
This concludes our induction proof.
\en
Consequently we have $ \tau (G) =G.$
Since $G\cmatrix {1 \cr 0 } = \cmatrix {\alpha \cr \beta },$
we see that 
$G= \cmatrix{ \alpha & \gamma \cr \beta &
  \lambda },$ with $ \beta = 2 \gamma.$
Since $ \det (G)= \pm 1,$ we obtain $ \beta ^2 \equiv \pm 2 \Mod\alpha.$
\EPf
\section{The harmonic knots ${\H(a,b,c)}$}\label{harmonic}
\begin{figure}[t!h]
\begin{center}
\begin{tabular}{ccc}
{\scalebox{.15}{\includegraphics{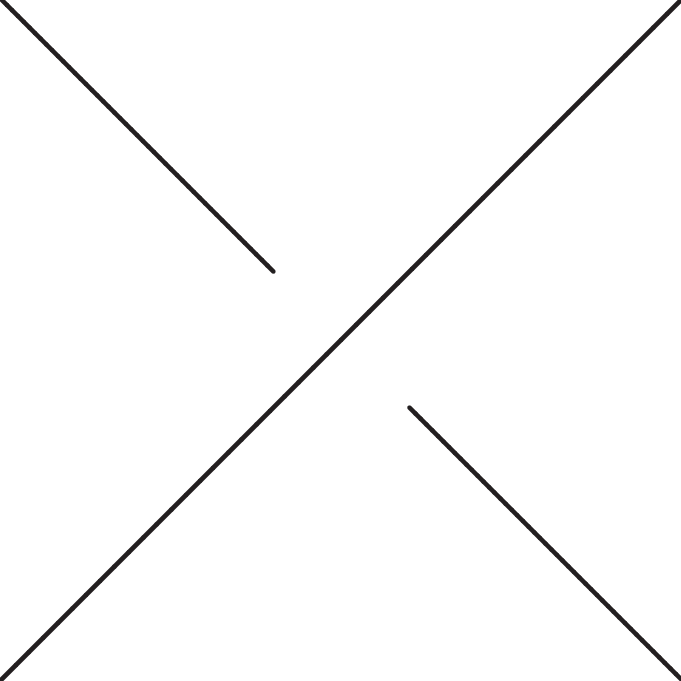}}}&\quad\quad &
{\scalebox{.15}{\includegraphics{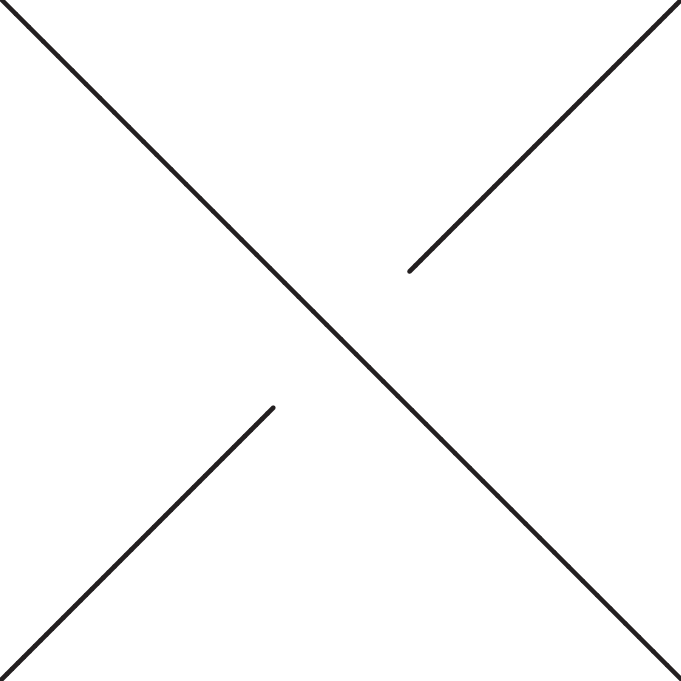}}}\\
A right twist & \quad & A left twist
\end{tabular}
\end{center}
\vspace{-12pt}
\caption{Right and left twists}
\label{fig:right-left}
\end{figure}

We shall first show some properties of the plane Chebyshev curves
$x=T_a(t), \, y=T_b(t)$. The following result is proved in
\cite{KP3}.
\begin{proposition}\label{prop:double-points}
Let $a $ and $b$ be coprime integers.  The $\frac 12 (a-1)(b-1)$
double points of the Chebyshev curve
$x= T_a(t), y= T_b(t)$ are obtained for the parameter pairs
$$
t= \cos \bigl( \Frac ka + \Frac hb \bigr) \pi,  \
s = \cos \bigl( \Frac ka - \Frac hb \bigr) \pi,
$$
where $h,k$ are positive integers such that
$ \Frac ka + \Frac hb < 1.$
\end{proposition}
Let us define a right twist and a left twist as in Figure \ref{fig:right-left};
this notion depends on the choice of the coordinate axes.
\pn
We shall need the following result:
\begin{lemma}[\cite{KP3,KPR}]\label{sign}
Let $\H(a,b,c)$ be a harmonic knot. \\
A crossing point $M$ of parameters $ (t,s)$
is a right twist if and only if
$$D(M) = \Bigl( z(t)-z(s) \Bigr) x'(t) y'(t) >0 \quad  . $$
\end{lemma}
From   Proposition \ref{prop:double-points}  and  Lemma \ref{sign}, we immediately deduce
\begin{corollary}\label{cprime}
Let $a,b,c$ be coprime integers. Suppose that the integer $c'$ satisfies
$ c' \equiv c  \Mod{2a} $ and $ c' \equiv -c \Mod{2b}.$
Then the knot $\H(a,b,c')$ is the mirror image of $\H( a,b,c).$
\end{corollary}
\Pf
At each crossing point  we have
$
T_{c'}(t) - T_{c'}(s) = - \Bigl( T_c(t) - T_c (s)  \Bigr).
$
\EPf

The next result is useful to reduce  the degree of a harmonic knot.
\begin{corollary} \label{reduc}
Let $a,b,c$ be coprime integers. Suppose that the integer $c$ is of
the form $c= \lambda a + \mu b$ with $\lambda, \mu >0$.
Then there exists $c'< c $
such that
$\H( a,b,c')$ is the mirror image of ${\H} (a,b,c). $
\end{corollary}
\Pf
Let $c'=\abs{\lambda a - \mu b}.$
The result follows immediately from Corollary \ref{cprime}.
\EPf
\pn
In \cite{KP4} we obtained the Schubert fractions of the harmonic knots
$ \H (3, b ,c )$, and their classification.
We shall follow the same strategy to study the harmonic knots $\H(4,b,c)$.
\subsection{The harmonic knots ${{\H(4,b,c)}}.$}
The following result characterizes the harmonic knots
 $\H(4,b,c).$
\begin{theorem}\label{th:h4}
Let $b, c$ be coprime odd integers such that $ b \not \equiv c  \Mod 4.$
The Schubert fraction of the knot $K= \H(4,b,c)$ is given by the
continued fraction
$$
\Frac \alpha\beta =
[ e_1, 2e_2, e_3, 2e_4, \ldots, e_{b-2}, 2e_{b-1} ],
{\rm where} \
e_i =  \sign { \sin (i -b) \theta} , \
\theta = \Frac{3b-c}{4b} \pi.$$
We have $\beta^2  \equiv \pm 2 \Mod \alpha.$
If $b< c < 3b,$
then the crossing number of $K$ is $N= ( 3b+c-2)/4.$
\end{theorem}
The proof will be given in section \ref{proofs}, p. \pageref{proofs}.
\pn
We are now able to classify  the harmonic knots of the form $\H(4,b,c)$.
\begin{theorem}\label{h4bc}
Let $b$ and $c$ be relatively prime odd integers, and let $K= \H(4,b,c).$
There is a unique pair
$(b',c')$
such that (up to mirroring)
$$K= \H(4,b',c'), \    b'<c'<3b', \  b' \not\equiv c' \Mod 4.$$
$K$ has a   Schubert fraction
$\Frac \alpha \beta $  such that
$ \beta^2 \equiv \pm 2 \Mod \alpha.$
Furthermore, there is an algorithm to find $(b', c'),$
and the crossing number of $K$ is $ N= (3b'+c'-2)/4. $
\end{theorem}
\Pf
First, let us  prove the uniqueness of this pair.
Let $K= \H(4,b,c)$ with $ b<c<3b, \ c \not\equiv b \Mod 4$.
By Theorem \ref{th:h4}, $K$ admits a Schubert fraction
$\Frac \alpha \beta $
such that $ \beta ^2 \equiv \pm 2 \Mod  \alpha, $ which implies that
$\alpha \ne 5 .$

Suppose that $\Frac \alpha{\beta'}$ is another Schubert fraction of
$K$ (or $ \overline{K}$) with $0< \beta'< \alpha$, $\beta'^2 \equiv \pm 2 \Mod \alpha.$
 We have $\beta\beta' \equiv \pm 1 \Mod {\alpha}$
so $\pm 4 \equiv 1 \Mod{\alpha}$.
Since $\alpha \ne 5,$ we see that $ \alpha = 3 ,$ and
then $  \beta =2 ,$ and $ \beta ' = 1$ is odd.

Consequently, there is a unique Schubert fraction $\Frac{\alpha}{\beta}$
of $K$ (or  $ \overline{K}$) such that
 $ 0< \beta < \alpha$,
$\beta^2 \equiv \pm 2 \Mod \alpha$ and $\beta$ even.
By Theorem \ref{th:h4}, the integer $b-1$ is  the
length of the  continued fraction expansion without three consecutive
sign changes of
$\Frac \alpha \beta = [e_1, 2e_2,  \ldots, e_{b-2}, 2e_{b-1} ]$.
Since we also have $3 b+c -2= 4 \, \cn(K),$
we deduce that the pair $(b,c)$ is uniquely determined.
\pn
Now, let us prove the existence of the pair $(b',c').$
Let $K= \H(4,b,c),  \   b<c$.
We  have only to  show that if the pair $(b,c)$
does not satisfy
the condition of the theorem, then it is possible to reduce it.
\pn
If $ c \equiv b \Mod 4, $ then $c= b + 4 \mu, \  \mu >0, $
and we can reduce the pair $ (b,c) $ by Corollary \ref{reduc}.
\pn
If $ c \not\equiv b \Mod 4 $ and $c>3b,$
then we have $c= 3b + 4 \mu, \  \mu >0,$ and we can reduce $(b,c)$ by Corollary \ref{reduc}.
\EPf
\begin{remark}
It follows  that
 the knots
$\H(4,b,c), \  4<b<c,  \  c \neq 4\lambda + \mu b, \  \lambda, \mu >0 $
are  different knots.
We also see that the only knot belonging to the two families
 $ \H(3,b,c)$ and $\H(4,b,c)$
is the trefoil $ \H(3, 4,5 )=  \overline  \H(4,3,5).$
\end{remark}
\begin{corollary}\label{h4special}
 The harmonic knot  $\H(4,2k-1,2k+1)$ is the two-bridge knot
of Conway form $ C(3,2, \ldots ,2)$ and crossing number $2k-1.$
\end{corollary}
\Pf
By Theorem \ref{th:h4}, the knot $\H_k=\H(4,2k-1,2k+1)$ has crossing number
$2k-1$ and Conway form
$
C (e_1,2e_2, \ldots, e_{2k-3},2e_{2k-2})$, where
$e_j =  \sign{ \sin ( j-b) \theta},$
 $\theta = \frac {\pi}2 (1- \frac{1}{2k-1}).$

Since the knots $C(a_1, \ldots , a_{2m})$ and $ C ( -a_{2m}, \ldots  ,
- a_1)$ are isotopic,
we deduce that $\H_k$ is isotopic to the knot
$ C ( 2 \eps_1, \eps_2, \ldots , 2 \eps_{2k-3}, \eps_{2k-2} )$
where
$ \eps_i= \sign { \sin i \theta } = (-1)^{\pent {i-1}2 }$.

We deduce that the rational number
$r_k= [2,1,-2,-1,\ldots, (-1)^{k-2}2, (-1)^{k-2}] $   (length $2k-2$) is
a Schubert fraction of $\H_k.$
We have $r_2=3,$ and $r_k = [2,1,-r_{k-1} ].$
Using the identity $ [2,1,-x]= [3, x-1],$ by  an easy induction we obtain
$ r_k= [3,2,\ldots, 2]. $
\EPf
\begin{example}[The twist knots]\label{twh4}
The  twist knots $\cT_n= C(n,2)$  are not harmonic knots $\H(4,b,c) $ for $n>3.$
They are not harmonic knots $ \H( 3,b,c)$ for $n >2.$
\end{example}
\Pf
The Schubert fractions of $\cT_n$
(or $ \overline{\cT_n} $)
with an even denominator are
$ \Frac {2n+1}2$, and  $ \Frac{2n+1}{-n} $ or $ \Frac {2n+1}{n+1}$ depending on
the parity of $n$.
The only such fractions satisfying $ \beta ^2 \equiv \pm 2 \Mod \alpha $ are
$\Frac 32$, $\Frac 74$ or $\Frac 94$.
The first two are the Schubert fractions of the trefoil and the $\overline{5}_2$ knot,
which are harmonic for $a=4.$
The case of $6_1 = S ( \Frac 9{4} )$ remains to be studied.
We have
$ \Frac 9{4} =  [ 1,2, -1,2, 1,-2, 1,2 ].$
Since this continued fraction expansion has two consecutive sign changes, by Theorems \ref{th:cf1212}
and \ref{th:h4}
we see that $6_1$ is not of the form
$\H(4,b,c).$
\EPf
\pn
But there also exist  infinitely many two-bridge knots whose Schubert fractions
$\Frac\alpha\beta$ satisfy $\beta^2\equiv -2 \Mod\alpha$  that are not
harmonic knots for $a=4$.
\begin{proposition}
The knots $S(n+\Frac{1}{2n})$ are not harmonic knots $\H(4,b,c)$ for $n>1$.
Their crossing number is $3n$ and their Schubert fractions $\Frac\alpha\beta = \Frac{2n^2+1}{2n}$
satisfy $\beta^2\equiv -2 \Mod \alpha$.
\end{proposition}
\Pf
We shall use the M\"{o}bius transformations
$$ F(x)= [1,2,x]= \Frac {3x+1}{2x+1},
C(x)= [1,2,-1,2,x]=x+2,
D(x)= [ 1,-2,1,2,x] = \Frac x {4x+1}.$$ 
We have $C^k(x) = 2k+x$ and $ D^k(x) = \Frac x{4kx+1}$, so $D^k (\infty)= \Frac{1}{4k}$.

If $n=2k$, we get $n+\Frac{1}{2n} = C^k D^k(\infty).$

If $n=2k+1$, we have
$n+\Frac{1}{2n} = n-1 + \Frac{2n+1}{2n} = C^k F D^k(\infty).$

These continued fractions are such that $ \beta ^2 \equiv -2 \Mod \alpha .$
Nevertheless, for $n>1$ these continued fractions  have two consecutive sign changes, and therefore they do not correspond to harmonic knots $\H(4,b,c)$.
\EPf
\subsection{Proof of theorem \ref{th:h4}}\label{proofs}
By Proposition \ref{prop:double-points} the parameters of the  crossing points of the plane projection of
$\H= \H (4,b,c)$ are obtained for the parameter pairs $ ( t, s) $ where
$$
t= \cos \bigl( \Frac k4 + \Frac hb \bigr) \pi,  \
s = \cos \bigl( \Frac k4 - \Frac hb \bigr) \pi,
$$
where $h,k$ are positive integers such that
$ \Frac k4 + \Frac hb < 1.$
If we define
$m=\abs{kb-4h}, \ m'=kb+ 4h,$
then we have
$ t=\cos \bigl( \Frac m'{4b} \pi \bigr), \ s= \cos \bigl( \Frac {m}{4b} \pi \bigr).$
We shall denote
$ \lambda = \Frac {3b-c}4$
(or $ c= 3b - 4 \lambda$), and $ \theta = \Frac {\lambda}b \pi.$

If $x,y$ are real numbers, then we shall write $x \sim y$ to mean that $ xy >0.$

We have to consider  two  cases.
\pn
{\bf The case $\mathbf{b = 4n+1}$.}
\pn
For $j=0, \ldots, n-1$,  let us consider the following crossing points
\bi
\item
$A_j$ corresponding to $ m=4j+1, \  m'= 2b-m,$ (or $ k=1, \ h= n-j $),
\item
$B_j$ corresponding to $ m= 4j+2, \  m'= 4b-m,$ (or $ k=2, \ h=2n-j$),
\item
$C_j$ corresponding to $ m= 4j+3, \ m' = 2b+m,$ (or $k=1, \ h=n+j+1$),
\item
$ D_j$ corresponding to $ m = 2b - 4(j+1), \  m' = 4b-m$ (or $ k=2, \ h=j+1$).
\ei
Then we have
\bi
\item
$ x( A_j)= \cos \bigl( \Frac {4j+1}b \pi \bigr), \  \  y(A_j) = (-1)^j \cos \Frac \pi 4 \neq 0, $
\item
$ x(B_j) = \cos \bigl( \Frac {4j+2}b \pi \bigr), \   \  y(B_j)=0,$
\item
$ x( C_j) = \cos \bigl( \Frac {4j+3}b \pi \bigr), \  \
 y(C_j) = (-1)^j \cos \Frac {3 \pi }4 \neq 0, $
\item
$x(D_j) = \cos \bigl( \Frac {4j+4}b \pi \bigr), \  \  y(D_j)=0.$
\ei
\psfrag{a0}{\small $A'_0$}\psfrag{b0}{\small $B_0$}\psfrag{c0}{\small $C_0$}\psfrag{d0}{\small $D_0$}%
\psfrag{an}{\small $A_{n-1}$}\psfrag{bn}{\small $B_{n-1}$}\psfrag{cn}{\small $C'_{n-1}$}\psfrag{dn}{\small $D_{n-1}$}%
\psfrag{aa0}{\small $A_0$}\psfrag{cc0}{\small $C'_0$}%
\psfrag{aan}{\small $A'_{n-1}$}\psfrag{ccn}{\small $C_{n-1}$}%
\begin{figure}[th]
\begin{center}
{\scalebox{.7}{\includegraphics{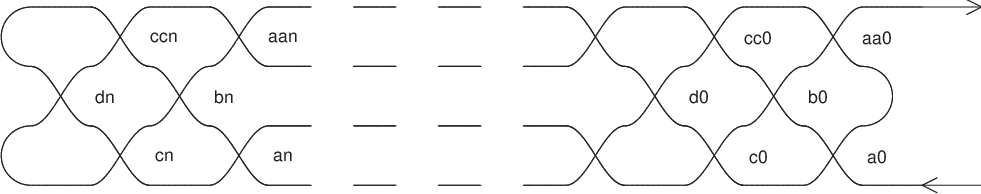}}}
\end{center}
\vspace{-12pt}
\caption{$\H(4,4n+1,c)$}
\label{dh4}
\end{figure}
Hence our $4n$ points satisfy
$$ x(A_0) > x(B_0) > x(C_0) > x(D_0) > \ldots >
x(A_{n-1}) > x( B_{n-1} ) > x( C_{n-1}) > x( D_{n-1}).$$
Let $A'_j$ (respectively $C'_j$) be the reflection of $A_j$
(respectively $C_j$) in the $x$-axis.
The crossings of our diagram are the points $ A_j, A'_j, B_j, C_j, C'_j, $ and $D_j.$
If $(t,s)$ is the parameter pair corresponding to
$A_j$ (respectively $C_j$), then $(-t,-s)$ is the parameter pair corresponding to $A'_j$ (respectively $C'_j$).
The   sign of  a crossing point  $M$ is
$s(M)= \sign{D(M)}$ if $y(M)=0,$ and
$s(M)=- \sign{D(M)}$ if $y(M) \neq 0.$
We have $s(A'_j) = s(A_j)$ and  $s(C'_j) = s(C_j).$

\begin{figure}[th]
\begin{center}
{\scalebox{.7}{\includegraphics{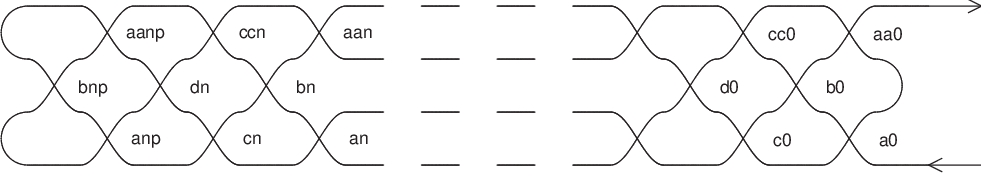}}}
\end{center}
\vspace{-12pt}
\caption{$\H(4,4n+3,c)$}
\label{dh42}
\end{figure}

A Conway  form of $\H$ is then (see section \ref{cf}, Figure \ref{conway4})
$$
C\Bigl( s (D_{n-1}), 2s(C_{n-1}), s(B_{n-1}),2s( A_{n-1}), \ldots, s (B_0), 2s (A_0) \Bigr).
$$
Using the identity $ T' _a ( \cos \tau ) = a \Frac { \sin a \tau }{ \sin \tau},$
we get
$ x'(t) y'(t)\sim \sin \bigl( \Frac mb \pi \bigr)  \sin \bigl( \Frac {m}4 \pi  \bigr).$
Consequently,
\bi
\item
For $A_j$ we have
$ x'(t) y'(t) \sim
\sin \bigl( \Frac {4j+1}b \pi \bigr) \sin \bigl( \Frac {4j+1}4 \pi \bigl)
\sim(-1)^j.$
\item
Similarly, for $B_j$, $ C_j$ and $D_j$ we obtain
$ x'(t) y'(t) \sim (-1)^j. $
\ei
On the other hand, at the crossing points we have
$$ z(t) - z(s) = 2 \sin \Bigl( \Frac c{8b} (m'-m) \pi \Bigr) \,
\sin  \Bigl( \Frac c{8b} (m+m') \pi  \Bigr).
$$
We obtain the  signs of our crossing points,
with $c = 3b-4 \lambda, \  \theta= \Frac {\lambda}b \pi.$
\bi
\item
For $A_j$ we get:
$ z(t )-z(s)= 2 \sin \Frac cb (n-j) \pi\, \sin   \Frac {c\pi}4.$\\
We have
$ \sin \Frac {c\pi}4 = \sin  \Frac { 12 n +3 - 4 \lambda}4  \pi =
(-1)^{n+ \lambda} \sin  \Frac {3 \pi }4 \sim (-1)^{n+ \lambda}$\\
and also
$
\begin{array}[t]{rcl}
\sin \bigl( \Frac cb (n-j) \pi  \bigr) &=&
\sin \Bigl(  \bigl( 3 - \Frac {4 \lambda}b \bigr) \bigl( n-j) \pi \Bigr)\\
&=& (-1)^{n+j} \sin \Bigl(  \Frac {4j-4n}b \lambda \pi  \Bigr) = (-1) ^{ n+j+ \lambda} \sin ( 4j+1) \theta
\end{array}$.\\
Consequently, the  sign of $A_j$ is
$ \ s(A_j) = - \sign{\sin ( 4 j+1) \theta }.$
\item
For $B_j$, we have:
$
\begin{array}[t]{rcl}
z(t) - z(s) &=& 2 \sin \bigl(  \Frac cb ( 2n-j) \pi \bigr)
\sin \Frac {c\pi}2 = - 2 \sin \bigl(  \Frac cb ( 2n-j) \pi \bigr)\\
&=& 2 \sin \Bigl( \bigl( 3 - \Frac{4 \lambda}b \bigr) \bigl( j-2n \bigr) \pi \Bigr) \\
&=& 2 (-1)^{j+1} \sin \Bigl( \Frac {\lambda}b ( 4j-8n) \pi \Bigr) = 2(-1)^{j+1} \sin (4j+2) \theta.
\end{array}
$
Therefore the  sign of $B_j$ is
$\ s (B_j) = - \sign{\sin (4j+2) \theta}.$
\item For $C_j$:
$ z(t) - z(s) = 2 \sin \Frac {c\pi}4
\sin \bigl(  \Frac cb (n+j+1) \pi  \bigr). $\\
We know that $ \sin \Frac {c \pi} 4 \sim (-1)^{n+ \lambda}$.
Let us compute the second factor:\\
$\sin \Bigl(  \bigl( 3 - \Frac {4 \lambda}b \bigr) \bigl(  n+j+1 \bigr) \pi\Bigr) =
(-1)^{n+j} \sin \Bigl( \Frac {\lambda}b \bigl( 4n+4j+4 \bigr) \pi \Bigr)$ \\
\hbox{} \hfill
$=(-1)^{n+j} \sin \Bigl( \Frac {\lambda}b ( b+4j+3) \pi \Bigr)=
( -1)^{n+j+ \lambda} \sin (4j+3)  \theta.$\\
Hence the  sign of $C_j$ is
$ s( C_j) = - \sign{\sin( 4j+3) \theta}.$
\item
For $D_j$:
$
\begin{array}[t]{rcl}
z(t) -z(s)&=& 2 \sin \bigl( \Frac cb (j+1) \pi \bigr) \sin \Frac {c\pi}2 \\
&=& -2 \sin \Bigl(  \bigl( 3- \Frac {4 \lambda }b   \bigr) \bigl( j+1 \bigr) \pi \Bigr) =
2(-1)^{j+1} \sin ( 4j+4) \theta .
\end{array}$.\\
We conclude that
$s (D_j)=
- \sign{\sin( 4j+4) \theta}.$
\ei
This completes the computation of our Conway  form of $\H$ in this first case.
\pn
{\bf The case $\mathbf{b = 4n+3}$.}\\
In this case  the diagram is different from the preceding one, see Figure \ref{dh42}.
As in the first case  the proof relies on carefully determining the sign of each crossing of the diagram.
The details are in \cite{KP5}.

\psfrag{aanp}{$A_n$}
\psfrag{anp}{$A'_n$}
\psfrag{bnp}{$B_n$}
In both cases  the  Conway  form of $\H(4,b,c)$ is
$
C (e_1,2e_2, \ldots, e_{b-2},2e_{b-1})
$
where $e_i =  \sign{\sin(i-b) \theta} .$
Consequently, we have  $ \beta ^2 \equiv \pm2 \Mod \alpha $ by
Proposition \ref{palin4}.
\pn
If $ b<c< 3b$ then we get
$ \lambda < \Frac b2,$ and   $ \theta < \Frac {\pi} 2.$
Consequently, there are no two consecutive sign changes in our sequence.
Moreover, the total number of sign changes is $ \lambda-1.$
We conclude by Proposition \ref{bireg} that the crossing number is
$ N= \Frac { 3 (b-1) }2 - (\lambda -1)= \Frac {3b+c-2}4.$
\EPf
\section{Some families with ${ a \ge 5}$} \label{special}
We will consider Chebyshev curves as trajectories in a rectangular
billiard (see \cite{KP3}).
\begin{lemma}\label{billiard}
Let  $\cC(t)$ be the plane curve parametrized by $ x(t)= T_a(t), \ y(t)=
T_b(t),$ and let $F$  be the function defined by
$ F(x) = \Frac 2{  \pi}  \arccos  x   -1 $.
The mapping $(x,y)\mapsto (X, Y)=(b\,F(x),a\,F(y))$ is a homeomorphism from
the square $I = (-1,1)^2$ onto the rectangle
$(-b, b) \times (-a , a)$.
The image of the curve  $\cC(t)$
is a ``billiard trajectory''. The slopes of its segments are $ \pm1,$
which means that they are parallel to one of the two  lines $ Y= \pm  X.$
\end{lemma}

\begin{figure}[th]
\begin{center}
\begin{tabular}{cccc}
{\scalebox{.8}{\includegraphics{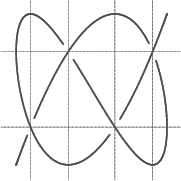}}}&
{\scalebox{.5}{\includegraphics{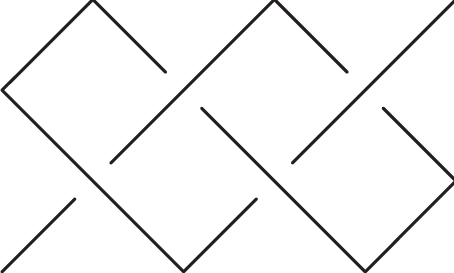}}}&
{\scalebox{.8}{\includegraphics{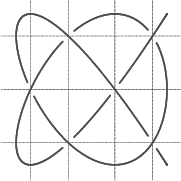}}}&
{\scalebox{.5}{\includegraphics{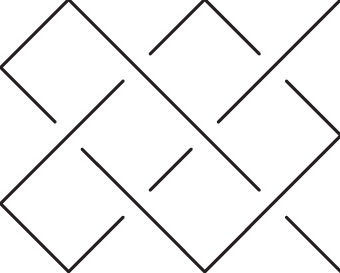}}}\\
\multicolumn{2}{c}{$\H(3,5,7)$}&\multicolumn{2}{c}{$\H(4,5,7)$}
\end{tabular}
\end{center}
\vspace{-12pt}
\caption{Billiard representations of $4_1$ and $5_2$}
\label{knbi}
\end{figure}
\subsection{The harmonic knots $\mathbf{ \H(2n-1, 2n,2n+1)}$}
Let us begin with some simple observations on the diagram of $K_n=
\H(2n-1,2n,2n+1)$.

\begin{figure}[t!h]
\begin{center}
\begin{tabular}{ccc}
{\scalebox{.7}{\includegraphics{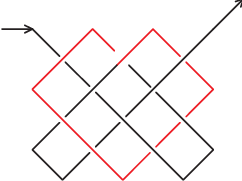}}}&
{\scalebox{.7}{\includegraphics{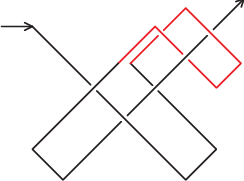}}}&
{\scalebox{.7}{\includegraphics{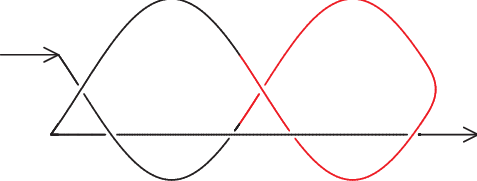}}}
\end{tabular}
\end{center}
\vspace{-12pt}\caption{An isotopy of $K_3$}
\label{k3}
\end{figure}

We have $z(t)=2t \, y(t) - x(t).$
Consequently, if $(t,s)$ is a parameter pair corresponding to a
crossing, we have: $ z(t)-z(s)= 2 (t-s)y(t)$.
This simple rule allows us to draw by hand the billiard picture
of the knot $K_n$ (see Figure \ref{kn}):
\begin{figure}[th]
\begin{center}
\begin{tabular}{ccc}
{\scalebox{.8}{\includegraphics{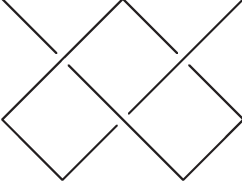}}}&
{\scalebox{.8}{\includegraphics{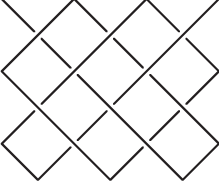}}}&
{\scalebox{.8}{\includegraphics{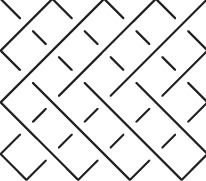}}}\\
$\H(3,4,5)$&$\H(5,6,7)$&$\H(7,8,9)$
\end{tabular}
\end{center}
\vspace{-12pt}\caption{The knots $K_n$ for $n=2,3,4$}
\label{kn}
\end{figure}
\pn
We can even deduce a simpler rule as follows.
\begin{lemma} Let $K= \H(a,b,c)$ with $b=a+1.$
Then the sign of a crossing point $M$ of parameters $(s,t)$ is
$  \sign {D(M)} = \sign{(z(t)-z(s))(t-s)}$.
\end{lemma}
\Pf
Let $(s,t)$ be the parameter pair of a crossing.
We have
$$
t=\cos \bigl(\Frac ka + \Frac hb \bigr) \pi  , \  s=\cos  \bigl( \Frac
ka - \Frac hb \bigr) \pi , \ 0 <  \Frac ka + \Frac hb < 1.
$$
An easy calculation shows that, when $b=a+1$ then
$$x'(t)y'(t)  \sim  - \sin ( \Frac ka \pi) \, \sin ( \Frac
hb \pi )  \sim t-s ,$$
which concludes the proof by using Lemma \ref{sign}.
\EPf
\begin{corollary} The sign of a crossing $M$ of $\H (  2n-1, 2n ,  2n+1)$ is
$\sign {D(M)} = \sign y .$
\end{corollary}
\begin{theorem} \label{surprise}
The knot $ \H (2n-1,2n, 2n+1)$ is isotopic to $ \H ( 4, 2n-1, 2n+1) $
if $n$ is odd, and to $ \H (4, 2n+1, 2n-1) $ if $n$ is even. Its
crossing number is $2n-1.$
\end{theorem}
\Pf
We shall use the billiard diagrams of harmonic knots defined in Lemma \ref{billiard}.
These diagrams are centered around the origin.
Our proof is   by induction on $n$.
We shall prove that $K_n$ is isotopic to the  two-bridge knot of Conway form
$ C(1,2,-1,-2, \ldots , (-1)^{n-2}, 2 (-1)^{n-2} ).$
\pn
For $n=2,$ the knot $\H(3,4,5)$ is the trefoil $K_2= C (1,2) =
\overline \H(4,3,5)$.
\pn
For $n=3,$ Figure \ref{k3}  shows that $K_3= C(1,2,-1,-2).$ It also gives an idea of our proof.
\pn
By induction, let us suppose that
$K_{n-1}= C( 1,2,-1,-2,\ldots, (-1)^{n-3}, (-1)^{n-3} \, 2 ).$
We shall consider that $K_n$ is  composed of two parts.

\begin{figure}[t!h]
\begin{center}
{\scalebox{.8}{\includegraphics{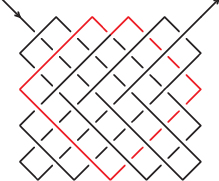}}} \quad
{\scalebox{.9}{\includegraphics{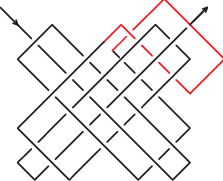}}} \quad
{\scalebox{.8}{\includegraphics{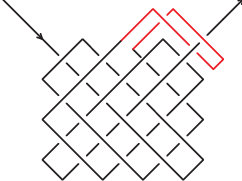}}}
\end{center}
\vspace{-12pt}\caption{Pulling the loop $L$  away from $K_n$.}
\label{kn1}
\end{figure}

The first part $L$ is a loop (the red loop of Figure \ref{kn1}) which
is symmetrical about the $y$-axis, and consists of the  points of  parameters
$ t \in I=( \pi (\frac 12 - \frac 1{2n-1}), \, \pi (\frac 12 + \frac
1{2n-1}) )$ It contains exactly $2(2n-3)$ crossing points, which are
the points of parameters
$$
t= \cos \tau, \  \tau = \Frac{\pi}2 + \Frac {k\pi}{2n(2n-1)} , \
 \abs k \le 2n-2, \  k \ne 0, \pm n.
$$
The other part $T_{n-1}$ consists of the points of parameters $t \in {\bf R} - I,$
it is a tangle over the rectangle $ (-2n, 2n) \times (-2n+1, 2n-1).$
\pn
When $n$ is odd,  the part of the loop  $L$ where $ t< \Frac \pi 2$
is over $T_{n-1}$,
and the other part of $L$  is under  $T_{n-1}.$
When $n$ is even, the first part of $L$ is under and the second part
of $L$  is over  $T_{n-1}$.
Consequently, it is possible to move the loop $L$ away from the
box containing  $T_{n-1}$ and we see that $K_n$ is obtained from $T_{n-1}$ by a
weaving process (see \cite[p. 50]{Ka}).

\begin{figure}[t!h]
\begin{center}
{\scalebox{.8}{\includegraphics{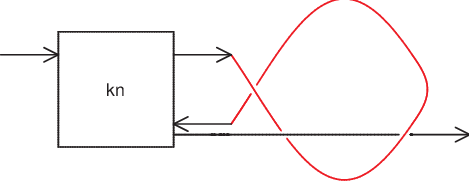}}} \quad
{\scalebox{.8}{\includegraphics{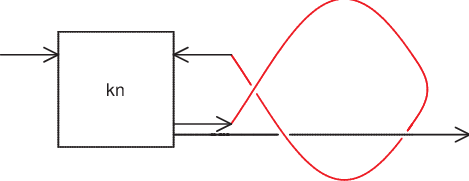}}}
\end{center}
\vspace{-12pt}\caption{The weaving process: $n$ odd (left), $n$ even (right)}
\label{weaving}
\end{figure}

Now let us look at the diagram of $T_{n-1}.$
It is clear (see Figure \ref{kn1}) that  the knot $K_{n-1} $ is the numerator
of the tangle $T_{n-1}$.
\pn
Consequently,  our weavings are illustrated in Figure \ref{weaving}.
\psfrag{kn}{$T_{n-1}$}
\pn
If $n$ is even, then using the induction hypothesis, we obtain the  Conway form
$ K_n= C(1,2, -1,-2, \ldots  , 1,2) $ of length $2n-2.$
If $n$ is odd, then  we obtain the Conway form
$ K_n = C(1,2, -1,-2, \ldots , -1,-2)$ of length $2n-2.$
This completes our induction proof.

\begin{figure}[t!h]
\begin{center}
\psfig{figure=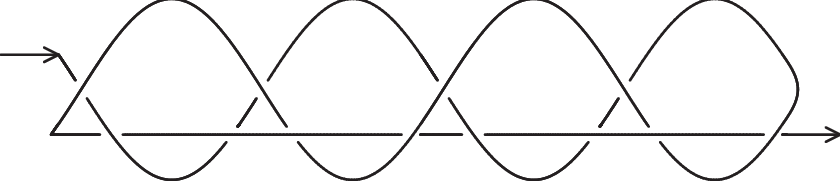,height=1.5cm,width=7cm} \quad
\psfig{figure=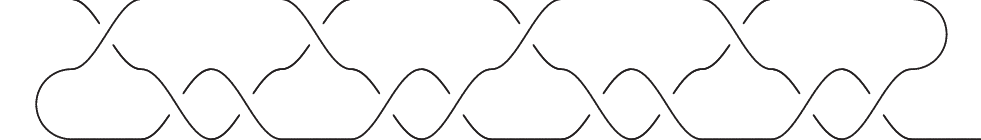,height=1.5cm,width=7cm}
\end{center}
\vspace{-12pt}\caption{The knot $K_5$ is a two-bridge knot}
\label{k5}
\end{figure}

By the proof of Corollary \ref{h4special}, we deduce that $K_n$ is isotopic
to $\H(4, 2n-1, 2n+1)$ if $n$ is odd, and to $\H(4, 2n+1, 2n-1)$ if
$n$ is even.
\EPf
\pn
The result of this inductive weaving process
is illustrated in Figure \ref{k5} for the knot $K_5.$
\subsection{The harmonic knots $ \mathbf{ \H (5, k, k+1).}$}\label{section:h5}
The bridge number of such a knot is at most three,
and one can verify
that the bridge number of
the knots $\H(5,5k+2,5k+3)$, $ \ 2 \leq k \leq 8$
is three.
This is the reason why the following result surprised us.
\begin{theorem} \label{h5}\\
The knot $ \H ( 5, 5n+1, 5n+2)$ is the two-bridge knot of Conway form $C ( 2n+1, 2n).$\\
The knot $ \H ( 5, 5n+3, 5n+4)$ is the two-bridge knot of Conway form $ C (2n+1, 2n+2).$\\
Besides $\H(5,6,7)= \H (4,5, 7)$ and $ \H(5,3,4),$
these knots are not of the form $ \H (a,b,c)$ with $a \le 4.$
\end{theorem}
The proof of this result is contained in \cite{KP5}. It is very similar to the preceding one.
\subsection{Some new findings on harmonic knots}
Thanks to the simplicity of our billiard diagrams, we can easily compute
the Alexander polynomials of our knots (see
\cite{Li}). On the other hand,  there is a list of the Alexander polynomials
of the first prime knots with 15 or fewer crossings  in \cite{KS}.

\begin{figure}[t!h]
\begin{center}
\begin{tabular}{ccccc}
{\scalebox{.6}{\includegraphics{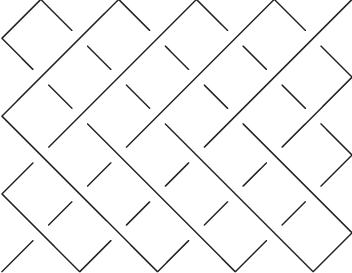}}} &\quad&
{\scalebox{.6}{\includegraphics{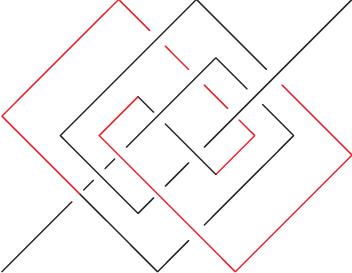}}} &\quad&
{\scalebox{.6}{\includegraphics{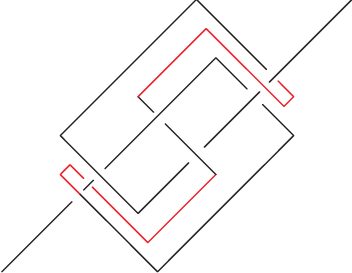}}}
\end{tabular}
\end{center}
\vspace{-12pt}\caption{The knot $\H(7,9,11),$ an unusual model of $8_{17}$}
\label{8-17}
\end{figure}

Using this list and some evident simplifications, we can
identify our knot.
We first give some specific examples, then an exhaustive list of
 knots $\H(a,b,c)$ having a  diagram
with 15 or fewer crossings.
\pn
{\bf  Harmonic knots are not necessarily prime.}\\
The knot $ \H( 5, 7, 11 )$ is not prime; it is the connected sum of
two figure-eight knots.
\begin{figure}[t!h]
\begin{center}
{\scalebox{.6}{\includegraphics{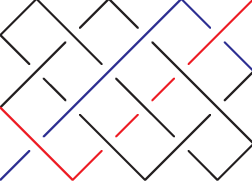}}} \quad
{\scalebox{.7}{\includegraphics{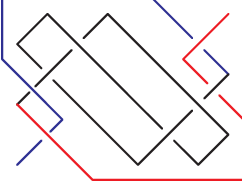}}}
\end{center}
\vspace{-12pt}\caption{The knot $\H ( 5,7,11)$ is composite}
\end{figure}
\pn
{\bf Harmonic knots may be nonreversible.}\\
We have identified the  knots of form
$\H (2n-1, 2n+1, 2n+3), \  n \le 5$, by computing their  Alexander
polynomials and  their crossing numbers.
We found two nonreversible harmonic knots, namely
$\H(7,9,11)  = 8_{17}$
 and
$\H(9,11,13)  = 10_{115}.$
\pn
Figure \ref{8-17} shows that $\H ( 7,9,11) = 8_{17}$ is symmetric
through the origin and therefore is strongly ($-$)amphicheiral. It is also
the first nonreversible knot (see  \cite[p.~30]{Cr}).
\pn
{\bf A table of harmonic knots  with $\mathbf{(a-1)(b-1) \le 30.}$}\\
Here, we provide a table giving the names (up to mirroring) of the  knots $\H(a,b,c)$
with diagrams having  $15$ or fewer crossings.
The knots are lexicographically ordered, and
by Corollary \ref{reduc} we choose $c$ such that
$c \not = \lambda a + \mu b, \, \lambda, \mu >0.$
 We have to identify 51   knots.

When $a=3$ or $a=4,$  $ \H(a,b,c)$ is a two-bridge
knot. The crossing number of such a knot is $\frac 13 (b+c)$, when
$a=3$ and $\frac 14 (3b+c-2)$ when $a=4$. Furthermore, its Schubert
fraction is computed using Theorem \ref{h4bc} or \cite[Theorem. 6.5]{KP4}.

When $a\geq 5,$ we compute the Alexander polynomial of the knot and
compare it with the tables. Sometimes (when starred) it is also
necessary to use their DT-notations and Knotscape (\cite{KS}).

\pn
\begin{center}
\begin{longtable}{||r|r|r||r|r|r||}
\hline
\multicolumn{6}{||c||}{{\bf Table of the first harmonic knots}}\\
\hline
& Fraction & Name &
& Fraction & Name\\
\hline
\H(3,4,5)&3&$3_1$&
\H(3,5,7)&5/2&$4_1$\\
\H(3,7,8)&5&$5_1$&
\H(3,7,11)&13/5&$6_3$\\
\H(3,8,13)&21/8&$7_7$&
\H(3,10,11)&7&$7_1$\\
\H(3,10,17)&55/21&$9_{31}$&
\H(3,11,13)&17/4&$8_3$\\
\H(3,11,16)&39/14&$9_{17}$&
\H(3,11,19)&89/34&$10_{45}$\\
\H(3,13,14)&9&$9_1$&
\H(3,13,17)&53/23&$10_{37}$\\
\H(3,13,20)&105/41&$11a175$&
\H(3,13,23)&233/89&$12a499$\\
\H(3,14,19)&77/34&$11a119^*$&
\H(3,14,25)&377/144&$13a1739$\\
\H(3,16,17)&11&$11a367$&
\H(3,16,23)&187/67&$13a2124^*$\\
\H(3,16,29)&987/377&$15a39533^*$&
\H(4,5,7)&7/2&$5_2$\\
\H(4,5,11)&11/3&$6_2$&
\H(4,7,9)&17/5&$7_5$\\
\H(4,7,13)&23/5&$8_7$&
\H(4,7,17)&41/11&$9_{20}$\\
\H(4,9,11)&41/12&$9_{18}$&
\H(4,9,19)&89/25&$11a180$\\
\H(4,9,23)&153/41&$12a541$&
\H(4,11,13)&99/29&$11a236$\\
\H(4,11,17)&113/31&$12a758$&
\H(4,11,21)&187/41&$13a2679^*$\\
\H(4,11,25)&329/87&$14a7552^*$&
\H(4,11,29)&571/153&$15a42637^*$\\
\H(5,6,7)&7/4&$5_2$&
\H(5,6,13)&&$10_{159}$\\
\H(5,6,19)&&$10_{116}$&
\H(5,7,8)&5/2&$4_1$\\
\H(5,7,9)&13/8&$6_3$&
\H(5,7,11)&&$4_1 \# 4_1$\\
\H(5,7,13)&&$12n356$&
\H(5,7,16)&&$12n798$\\
\H(5,7,18)&&$12n321$&
\H(5,7,23)&&$12a960$\\
\H(5,8,9)&13/4&$7_3$&
\H(5,8,11)&21/13&$7_7$\\
\H(5,8,17)&&$14n22712^*$&
\H(5,8,19)&&$14n26120^*$\\
\H(5,8,27)&&$14a19221^*$&
\H(6,7,11)&&$10_{134}$\\
\H(6,7,17)&&$15n42918^*$&
\H(6,7,23)&&$15n165258^*$\\
\H(6,7,29)&&$15a81117$&&&\\
\hline
\end{longtable}
\end{center}
\pn
We have observed that for some integers $b, \,  k$ and $h,$
$ \H(b-k,b,b+k) = \H(b-h,b,b+h)$.
It is the case for
$\H(5, 11, 17) = \H(9, 11, 13)$,
$\H(3, 11, 19)=\H(7, 11, 15)$
and many others. It would be interesting to explain this phenomenon.
\pn

\vfill
\pn\goodbreak
\hrule width 5cm height 2pt
\pn\nobreak
Pierre-Vincent Koseleff, \\
Universit{\'e} Pierre et Marie Curie (UPMC-Paris 6) \&
\'Equipe INRIA Ouragan \& Institut de Math{\'e}matiques de Jussieu (UMR-CNRS 7586)\\
e-mail: {\tt koseleff@math.jussieu.fr}
\pn
Daniel Pecker, \\
Universit{\'e} Pierre et Marie Curie (UPMC-Paris 6) \&
Institut de Math{\'e}matiques de Jussieu (UMR-CNRS 7586)\\
e-mail: {\tt pecker@math.jussieu.fr}
\end{document}